\documentclass[runningheads]{llncs}
\usepackage[fleqn]{amsmath}
\usepackage{amssymb}
\usepackage{graphicx}

\begin{document}

\mainmatter

\title{An Infinite Familiy of Quadratic Quadrinomial APN Functions}

\author{Carl Bracken \thanks{Research supported by
Irish Research Council for Science, Engineering and Technology
Postdoctoral Fellowship} 
\and
{Eimear Byrne \and Nadya Markin \and Gary McGuire} \thanks{Research supported by the Claude
Shannon Institute, Science Foundation Ireland Grant 06/MI/006}}

\institute{School of Mathematical Sciences\\
University College Dublin\\
Ireland\\
\email{carlbracken@yahoo.com, ebyrne@ucd.ie, nadyamarkin@gmail.com, garymcguire@ucd.ie}}


\maketitle

\begin{abstract}
\noindent We present an infinite familiy of APN functions on
$GF(2^{3k})$ with $(k,3)=1$.
\end{abstract}


\bigskip


\bigskip


\section{Introduction}

Let $L=GF(2^n)$ for some positive integer $n$.
A function $f : L \longrightarrow L$ is said to be \emph{almost perfect nonlinear} (APN) on $L$
if the number of solutions in $L$ of the equation
$$f(x+q)+f(x)=p$$
is at most 2, for all $p,q\in L$, $q \not=0$. Equivalently, $f$ is APN if the set
$\{f(x+q)+f(x): x \in L \}$ has size $2^{n-1}$ for each $q \in L^*$. 
Clearly, as $L$ has characteristic 2, the number of solutions
to the above equation must be an even number for any function $f$ on $L$.

APN functions were introduced in \cite{N} by Nyberg, who defined them as the mappings with
highest resistance to differential cryptanalysis. In other words, APN functions are those for which
the plaintext difference $x+y$ yields the ciphertext difference $f(x)+f(y)$ with probability $1/2^n$.
Since Nyberg's characterization, many papers have been written on APN functions, although not many
different families of such functions are known.

Two functions $f,g : L \longrightarrow L$
are called {\em extended affine} (EA) equivalent
if there exist affine permutations $A_1,A_2$ and an affine map $A$ such that
$g=A_1\circ f \circ A_2 + A$.

Until recently, all known APN functions were
EA equivalent to one of a short list of monomial functions,
namely the Gold, Kasami-Welch, inverse, Welch, Niho and Dobbertin functions.
For some time it was conjectured that this list was the complete
list of APN functions up to EA equivalence.

A more general notion of equivalence has been suggested in \cite{CCZ},
which is referred to as Carlet-Charpin-Zinoviev (CCZ) equivalence.
Two functions are called CCZ equivalent if the graph
of one can be obtained from the graph of the other by an affine permutation
of the product space. EA equivalence is a special case of CCZ equivalence.

We say that $f:L \longrightarrow L$ is differentially $m-$uniform if the polynomial
$f(x+q)+f(x)+p$ has at most $m$ zeroes in $L$. Then $f$ is APN on $L$ if and only if it is 
differentially 2-uniform on $L$. 

Differential uniformity and resistance to linear and differential attacks are 
invariants of CCZ equivalence, and as opposed to EA equivalence, any
permutation is CCZ equivalent to its inverse.

In \cite{BCP}, Proposition 3, the authors express necessary and sufficient
conditions for EA equivalence of functions in terms of CCZ equivalence
and use this to construct several examples of APN functions that are
CCZ equivalent to the Gold functions, but not EA equivalent to any monomial function.
This showed that the original conjecture is false.
The new question was whether all APN functions
are CCZ equivalent to one on the list.

In 2006 a sporadic example of a binomial APN function
that is not CCZ equivalent to any power mapping was given in \cite{EKP}.
A family of APN binomials on fields $\mathbb{F}_{2^n}$, where $n$ is divisible
by $3$ but not $9$, was presented in \cite{BCFL}.
In \cite{BCL} these have been shown to be
EA inequivalent to any monomial function, and CCZ inequivalent to the Gold
or Kasami-Welch functions.
For the case $n=6$, in \cite{Dillon} Dillon presented a list of CCZ inequivalent APN functions on $GF(2^n)$, found by computer search. In general, establishing CCZ equivalence of arbitrary functions is extremely difficult. There are, however, a number of invariants of CCZ equivalence that can be useful in the classification of functions. A nice link with coding theory is that a pair of functions $f$ and $g$ on $L$ are CCZ equivalent on $L$ if and only if the matrices 
$$H_f = \left[ \begin{array}{ccc}  {\bf x}_1 & \cdots {\bf x}_{2^n}\\ & \\  f({\bf x}_1) & \cdots  f({\bf x}_{2^n})\end{array} \right] ,\:\: H_g= \left[ \begin{array}{ccc}  {\bf x}_1 & \cdots {\bf x}_{2^n}\\ & \\  g({\bf x}_1) & \cdots  g({\bf x}_{2^n})\end{array} \right] $$
are parity check matrices for codes whose extended codes are equivalent over $GF(2)$, where ${\bf x}_i,f({\bf x}_i)$ and $g({\bf x}_i)$ are expressions of $x_i,f(x_i)$ and $g(x_i)$ respectively as binary vectors of length $n$ in $L$ viewed as a $GF(2)$ vector space and $L = \{x_1,...,x_{2^n}\}$. 

Motivated by these works, in this paper we introduce a new family of APN functions
on fields of order $2^{3k}$ where $k$ is not divisible by 3. The family of polynomials has the form
\begin{equation}\label{eqf}
F(x)=u^{2^k}x^{2^{-k}+2^{k+s}}+ux^{2^{s}+1}+vx^{2^{-k}+1}+wu^{{2^k}+1}x^{2^{k+s}+2^s}
\end{equation}
with certain constraints on the integers $s,k$ and on $u,v,w \in GF(2^{3k})$ (see Theorem \ref{thmain}).

In the next section we show that the polynomials of type (\ref{eqf}) are indeed APN on $GF(2^{3k})$.
Using code equivalence, in Section 3 we show that for $n=6$ the functions are CCZ inequivalent to any known power functions and are equivalent to one of the trinomials listed by Dillon in \cite{Dillon}. 

\section{New APN functions}

The following theorem will show that we can obtain quadratic
quadrinomial APN functions on $GF(2^n)$
 whenever $n$ is divisible by three but not nine. A quadratic monomial is one of the form
 $x^{2^i+2^j}$  for some integers $i$ and $j$. Observe that if $f(x)=x^{2^i+2^j}$, then
 $$f(x+q)+f(x)+f(q) = x^{2^i}q^{2^j}+x^{2^j}q^{2^i}$$
 is a linear function in $x$, whose kernel has the same size as any of its translates, such as the solution set of $f(x)+f(x+q)=p$ in $L$, for any $p \in L$. 
  
Note that because of this property, proving whether or not a quadratic polynomial is APN is more tangible than one that is not quadratic. For this reason, all of the recently discovered families of APN functions inequivalent to power mappings have been quadratic.

 We will show that our polynomial $F(x)$ is APN by computing the size of the kernel of the corresponding linear map $$F(x+q)+F(x)+F(q).$$

\begin{theorem}\label{thmain}
Let $s$ and $k$ be positive integers with $k+s$ divisible by three and
$(s,3k)=(3,k)=1$. Let $u$ be a primitive element of ${GF(2^{3k})}$ and let $v,w \in
GF(2^{k})$ with $v \neq w^{-1}$. Then the function 
$$F(x)=u^{2^k}x^{2^{-k}+2^{k+s}}+ux^{2^{s}+1}+vx^{2^{-k}+1}+wu^{{2^k}+1}x^{2^{k+s}+2^s}$$
is APN over $GF(2^{3k})$.
\end{theorem}

Proof:

We show that for every $p$ and $q$ (with $q \neq 0$) in $GF(2^{3k})$ the equation
$$F(x)+F(x+q)=p$$
has at most two solutions by counting the number of solutions to the equation
$$F(x)+F(x+q)+F(q)=0.$$
This gives
\begin{eqnarray*}
F(x)+F(x+q)+F(q) & = & u^{2^k}(x^{2^{k+s}}q^{2^{-k}}+q^{2^{k+s}}x^{2^{-k}})+u(x^{2^{s}}q+q^{2^{s}}x)\\
                                & + & v(x^{2^{-k}}q+q^{2^{-k}}x)+wu^{2^k+1}(x^{2^s}q^{2^{k+s}}+q^{2^s}x^{2^{k+s}})=0.
\end{eqnarray*}
Replace $x$ with $xq$ to obtain
\begin{eqnarray*}
u^{2^k}q^{2^{-k}+2^{k+s}}(x^{2^{k+s}}+x^{2^{-k}})+uq^{2^{s}+1}(x^{2^{s}}+x)
+  vq^{2^{-k}+1}(x^{2^{-k}}+x) \\ + wu^{2^k+1}q^{2^{k+s}+{2^s}}(x^{2^s}+x^{2^{k+s}})=0,
\end{eqnarray*}
and collect terms in $x$ to get
\begin{eqnarray*}
\Delta(x)  :=  (vq^{2^{-k}+1} +uq^{2^s+1})x+(vq^{2^{-k}+1}+u^{2^k}q^{2^{-k}+2^{k+s}})x^{-k} \:\:\:\:\:\:\:\\ 
                 +   (wu^{2^{k}+1}q^{2^{k+s}+2^s}  + uq^{2^s+1})x^{2^s} +(wu^{2^k+1}q^{2^{k+s}+2^s}+u^{2^k}q^{2^{-k}+2^{k+s}})x^{k+s}=0.
\end{eqnarray*}
We write 
$$\Delta(x)=Ax+Bx^{2^{-k}}+Cx^{2^s}+Dx^{2^{k+s}}$$ 
where 
\begin{eqnarray*}
A &= & vq^{2^{-k}+1}+uq^{2^s+1},\:\:\:\: 
B=vq^{2^{-k}+1}+u^{2^k}q^{2^{-k}+2^{k+s}},  \\
C & = & wu^{2^{k}+1}q^{2^{k+s}+2^s}+uq^{2^s+1},\:\:\:\:
D=wu^{2^k+1}q^{2^{k+s}+2^s}+u^{2^k}q^{2^{-k}+2^{k+s}}.
\end{eqnarray*} 
Clearly $0$ is a root of $\Delta(x)$. Moreover $\Delta(1) = A+B+C+D=0$. 
If we show that $\Delta(x)=0$ permits only 0 and 1 as solutions for $x$ then
we will have proved that $F(x)$ is APN on $GF(2^{3k})$. 
First we demonstrate that none of $A, B, C $ or $D$
vanish for any $q \in GF(2^{3k})^*.$ If $A=0$ we have
$u=vq^{2^{-k}-2^s}$ which implies $u^{2^k}=vq^{1-2^{k+s}}$. 
By hypothesis, $k+s$ is divisible by three so that $1-2^{k+s}$ is divisible by seven, and hence $q^{1-2^{k+s}}$ is a $7$th power
in $GF(2^{3k})$. Since $3$ does not divide $k$, $7$ does not divide $2^k-1$, 
so the map $x \mapsto x^7$ is a permutation on $GF(2^k)$. Then $v \in GF(2^k)$
can be expressed as a 7th power. This means that $u^{2^k}$ and hence $u$ is a 7th power in 
$GF(2^{3k})$. This gives a contradiction as
seven is a divisor of $2^{3k}-1$ and we chose u to be primitive in $GF(2^{3k})$. We deduce that 
$A \neq 0$. 
Similar arguments show that $B, C$ and $D$ are all nonzero. 

Next we define the linearized polynomial: 
$$L_\theta(T):=T+\theta T^{2^{k}}+{\theta}^{2^{k}+1}T^{2^{-k}}.$$
When $T=\theta x+ x^{2^{-k}}$ and $\theta$ is a $({2^{k}-1})$-th power,
a routine calculation verifies that $L_\theta(T)=0$ for all $x \in
GF({2^{3k}}).$ 
Observe that 
$$\frac{A}{B}=
\frac{vq^{2^{-k}+1}+uq^{2^s+1}}{vq^{2^{-k}+1}+u^{2^k}q^{2^{-k}+2^{k+s}}}
=\frac{v+uq^{2^s-2^{-k}}}{v+u^{2^k}q^{2^{k+s}-1}}=
({v+uq^{2^s-2^{-k}}})^{1-2^k},$$ 
which gives
\begin{equation}\label{eqL}
L_\frac{A}{B}\left(\frac{A}{B}x+x^{2^{-k}}\right)=0.
\end{equation}
Now $$ \frac{\Delta(x)}{B} = (\frac{A}{B}x + x^{2^{-k}}) +(\frac{C}{B}x^{2^s} + \frac{D}{B}x^{k+s}) =0 . $$
Applying this to Equation \ref{eqL} gives
$$L_\frac{A}{B}\left(\frac{\Delta(x)}{B}\right)=L_\frac{A}{B}\left(\frac{C}{B}x^{2^s}+\frac{D}{B}x^{2^{k+s}}\right)=0.$$
We compute this as
$$(B^{2^{-k}+{2^k}}C+D^{2^{-k}}A^{2^{k}+1})x^{2^s}+(B^{2^{-k}+{2^k}}D+B^{2^{-k}}AC^{2^k})x^{2^{k+s}}$$
$$+(B^{2^{-k}}AD^{2^k}+A^{2^k+1}C^{2^{-k}})x^{2^{-k+s}}=0.$$
We substitute in the values of $A, B, C, $ and $D$ and after
simplification we obtain the following
$$(vw+1)uq^{2^k+1+2^s}(vq^{2^{-k}}+uq^{2^s})(u^{2^k}q^{2^{k+s}+2^k}+u^{2^{-k}}q^{2^{-k+s}+1})x^{2^s}$$
$$+(vw+1)u^{2^k}q^{2^k+1+2^{k+s}}(vq^{2^{-k}}+uq^{2^s})(uq^{2^{k}+2^s}+u^{2^{-k}}q^{2^{-k+s}+2^{-k}})x^{2^{k+s}}$$
$$+(vw+1)u^{2^{-k}}q^{2^k+1+2^{-k+s}}(vq^{2^{-k}}+uq^{2^s})(u^{2^k}q^{2^{k+s}+2^{-k}}+uq^{2^s+1})x^{2^{-k+s}}=0.$$
As we chose $v$ and $w$ such that $v \neq w^{-1}$ and as $A \neq 0$
we can divide the equation by
$(vw+1)q^{2^k+1}(vq^{2^{-k}}+uq^{2^s})u^{2^{-k}+1}q^{2^{-k+s}+2^s+1}$
and take the expression to the $2^{-s}-th$ power to obtain
$$(1+a^{-2^{k-s}})x+(a^{2^{-s}}+a^{-2^{k-s}})x^{k}+(1+a^{2^{-s}})x^{2^{-k}}=0,  \ \ \ \ \  (1)$$
where $a=u^{2^{k}-1}q^{2^{-k}+2^{k+s}-2^{s}-1}.$ Now we consider
$L_\frac{C}{D}(\frac{\Delta(x)}{D})=0$. We know
$L_\frac{C}{D}(x^{2^s}+\frac{C}{D}x^{2^{k+s}})=0$, as
$$\frac{C}{D}=
\frac{wu^{2^{k}+1}q^{2^{k+s}+2^s}+uq^{2^s+1}}{wu^{2^k+1}q^{2^{k+s}+2^s}+u^{2^k}q^{2^{-k}+2^{k+s}}}
= ({w+u^{-1}q^{2^{-k}-2^s}})^{2^k-1}.$$ This implies
$L_\frac{C}{D}(\frac{A}{D}x+\frac{B}{D}x^{2^{-k}})=0$, which we
compute as
$$(C^{2^{-k}+2^k}A+C^{2^{-k}}DB^{2^k})x+(C^{2^{-k}}DA^{2^k}+D^{2^k+1}B^{2^{-k}})x^k$$
$$+(C^{2^{-k}+2^k}B+D^{2^k+1}A^{-k})x^{2^{-k}}=0.$$ A similar
computation to the one used above will yield
$$(1+a^{-2^{-k}})x+(1+a)x^{2^k}+(a+a^{-2^{-k}})x^{2^{-k}}=0. \ \ \ \ \  (2)$$
Now we combine equations $(1)$ and $(2)$ such that the terms in
$x^{2^{-k}}$ cancel. This will give
$$((1+a^{-2^{k-s}})(a+a^{-2^{-k}})+(1+a^{-2^{-k}})(1+a^{-s}))x+$$
$$((a^{2^{-s}}+a^{-2^{k-s}})(a+a^{-2^{-k}})+(1+a)(1+a^{-s}))x^{2^k}=0$$
which is the same as
$$((1+a^{-2^{k-s}})(a+a^{-2^{-k}})+(1+a^{-2^{-k}})(1+a^{2^{-s}}))(x+x^{2^k})=0.$$
If we show that
$(1+a^{-2^{k-s}})(a+a^{-2^{-k}})+(1+a^{-2^{-k}})(1+a^{2^{-s}}) \neq
0$ for all possible values of
 $a$ then we could conclude that $x \in GF(2^k)$.
To this end we consider the expression
$$(1+a^{-2^{k-s}})(a+a^{-2^{-k}})=(1+a^{-2^{-k}})(1+a^{2^{-s}}).$$
Rearranging we obtain
$$a=\frac{{(1+a^{-1})}^{2^{-k}}}{{(1+a^{-1})}^{2^{k-s}}}\frac{{(1+a)}^{2^{-s}}}{{(1+a)}^{2^k}}.$$
This implies $a$ is a $(2^{k+s}-1)$-th power which in turn implies
that it is a seventh power. As
$a=u^{2^{k}-1}q^{2^{-k}+2^{k+s}-2^{s}-1}=u^{2^{k}-1}q^{(2^{k+s}-1)(1-2^{-k})}$
we see that if $a$ is a seventh power then so is $u^{2^{k}-1}$ but
this is not possible as $k$ is not divisible by three and $u$ is
primitive. We can now state that all solutions to $\Delta(x)=0$ are
in $GF(2^k)$. Applying this to our original expression for
$\Delta(x)$ gives
$$(uq^{2^s+1}+u^{2^k}q^{2^{-k}+2^{k+s}})(x+x^{2^{s}})=0.$$
If $uq^{2^s+1}+u^{2^k}q^{2^{-k}+2^{k+s}}=0$ then $a=1$, but $1$ is a
seventh power, hence $(x+x^{2^{s}})=0$ which implies $x=0$ or $1$ as
$s$ is relatively prime to $3k$.

\section{The Case $n=6$}

For the case $n=6$ the polynomials introduced here takes one of the following forms:
$$u x^3 +vu^5 x^{10} +v x^{17} +u^4 x^{24} $$
$$u x^3 +v x^{17} +u^4 x^{24} $$
$$ u x^3 +vu^5 x^{10} +u^4 x^{24}$$
$$u x^3 + u^4 x^{24}, $$
for some primitive element $u \in GF(2^{6})$ and $v \in GF(4)$.
In the first 3 cases, the polynomials are CCZ equivalent to
$$x^3+x^{10}+ u x^{24} ,$$  which appears in Dillon's list, and in the last instance the polynomial
is CCZ equivalent to $x^3$.

\pagebreak

\end{document}